\documentclass{pnastwo}
\usepackage{amssymb,amsfonts,amsmath}
\usepackage{epsfig}
\usepackage[arrow]{xy}

\usepackage{latexsym,cite,epsf}
\usepackage{graphicx,color}

\def\Acal{\mathcal{A}} 
\def\Fcal{\mathcal{F}} 
\def\Ucal{\mathcal{U}} 
\def\CC{\mathbb{C}} 
\def\ZZ{\mathbb{Z}} 
\def\zz{\mathbf{z}}

\def\SL{\operatorname{SL}}
\def\xx{\mathbf{x}}

\issuedate{2014}
\volume{}
\issuenumber{}

\begin{document}

\title{Webs on surfaces, rings of invariants, and clusters
}

\author{Sergey Fomin\affil{1}{Department of Mathematics, 
University of Michigan, Ann Arbor, MI 48109}
\and
Pavlo Pylyavskyy\affil{2}{Department of Mathematics, University of Minnesota,
Minneapolis, MN 55414}
}

\contributor{Submitted to Proceedings of the National Academy of Sciences
of the United States of America}



\maketitle

\begin{article}

\begin{abstract}
We construct and study cluster algebra structures in rings of 
invariants of the special linear group action on collections of
three-dimensional vectors, co\-vectors, and matrices. 
The construction uses Kuperberg's calculus of webs on
marked surfaces with boundary. 
\end{abstract}


\dropcap{T}he special linear group $\SL(V)$ 
of a finite-dimensional complex
vector space $V$ en\-dow\-ed with a volume form naturally acts 
on (the coordinate ring of the space of) 
collections of vectors in~$V$, covectors in~$V^*$, 
and operators in~$SL(V)$. 
The ring of invariants for this action
conjecturally carries a \emph{cluster algebra} structure,
and typically many of them. 
In this paper, we focus on the case when $V$ is three-dimensional. 
In this case, we describe a class of cluster structures on~$V$
defined in terms of combinatorial topology of 
bordered surfaces with marked points on the boundary. 
The key role is played by a construction of a family of invariants
associated with web diagrams, in the sense of G.~Kuperberg. 

The paper is organized as follows. We first review \hbox{basic}
background on invariant rings and cluster algebras. We then
introduce the class of marked surfaces used in our construction,
and describe the relevant variation of Kuperberg's dia\-grammatic
calculus. 
Next follows a technical description of 
cluster structures in the rings of invariants under consideration.
We then state our main results: 
the assertion~that~our construction produces a well
defined cluster algebra structure; 
the invariance of this structure under the choices that the
construction depends on; 
and its basic functoriality properties. 
We then discuss a few key examples, 
formulate several conjectures,
and point out a connection to the fundamental work by V.~Fock and
A.~Goncharov. 

\vspace{-.1in}

\section{1. Rings of invariants}

Let $V\cong\CC^k$ be a vector space endowed with a volume form,
which induces a dual volume form on~$V^*$.
The special linear group $\SL(V)$ acts 
on both~$V$ and~$V^*$;
to make the latter a left action, one sets $(gu^*)(v)=u^*(g^{-1}(v))$, 
for $v\in V$, $u^*\in V^*$, and $g\in\SL(V)$. 
The group $\SL(V)$ also acts on itself, 
via conjugation. 
In this paper, we focus our attention 
on the ring 
\[
R_{a,b,c}(V)=\CC[(V^*)^a\times V^b \times (\SL(V))^c]^{\SL(V)}
\]
of $\SL(V)$-invariant polynomials on $(V^*)^a\times V^b \times (\SL(V))^c$. 
%
The closely related $\SL(V)$ action on
$\CC[(V^*)^a\times V^b \times\operatorname{End}(V)]$ was studied by
  Procesi~\cite{procesi}. 

\begin{theorem}[{\rm cf.\ \cite[Theorem 12.1]{procesi}}]
\label{thm:procesi}
 The ring of invariants $R_{a,b,c}(V)$ is generated by: 
\begin{itemize}
 \item the traces $tr(X_{i_1} \dotsc X_{i_r})$ of arbitrary
 (non-commutative) monomials in the $c$ matrices in $\SL(V)$;
 \item the pairings $\langle v_i, M w_j \rangle$, where $v_i$ is a
   vector, $w_j$ is a covector and $M$ is any monomial as before;
 \item volume forms $\langle M_1 v_{i_1}, \ldots, M_n v_{i_n}\rangle$,
 where $M_i$-s are monomials as before and $v_i$-s are vectors; 
 \item volume forms $\langle M_1 w_{i_1}, \ldots, M_n w_{i_n}\rangle$,
 where $M_i$-s are monomials as before and $w_i$-s are covectors. 
\end{itemize}
\end{theorem}

While we know~\cite{procesi} that a finite subset of the invariants
listed above generates the
ring $R_{a,b,c}(V)$, the minimal size of such subset is not
known, except for some special cases. 
An even harder problem is to give an explicit version of
the ``second fundamental theorem,'' 
describing the ideal of relations satisfied by a minimal
generating set. 

The case $c=0$ of Theorem ~\ref{thm:procesi}
(invariants of vectors and covectors) goes back to
H.~Weyl~\cite{weyl}. 
In the case $a=c=0$, 
one recovers a ``Pl\"ucker ring,''
the homogeneous coordinate ring of a Grassmannian
$\operatorname{Gr}_b(V)$ with respect to its Pl\"ucker embedding. 
Pl\"ucker rings are among the most important and thoroughly
studied examples of cluster algebras, see
\cite{scott, gsv-book}. 





\vspace{-.12in}

\section{2. Cluster algebras}

Cluster algebras \cite{ca1,ca4} are commutative rings endowed with a 
combinatorial structure of a particular kind.
For the cluster algebras studied in this paper,
the defining combinatorial data are encoded in a
\emph{quiver}~$Q$, a finite oriented loopless graph
with no oriented 2-cycles. 
Some vertices of $Q$ are designated as
\emph{mutable}; 
the remaining ones are called \emph{frozen}. 

Let $z$ be a mutable vertex in a quiver~$Q$.
The \emph{quiver mutation} $\mu_z$ transforms $Q$ into the new
quiver~$Q'=\mu_z(Q)$ defined as follows.
First, for each pair of directed edges $x\to z\to y$
passing through~$z$, we introduce a new edge $x\to y$ (unless both
$x$ and~$y$ are frozen, in which case do nothing). 
Next, we reverse the direction of all edges incident
to~$z$. 
We then remove all oriented 2-cycles to obtain~$Q'$.

The combinatorial dynamics of quiver mutations drives the algebraic
dynamics of \emph{seed mutations}. 
Let~$\Fcal$ be a field containing~$\CC$. 
A \emph{seed} in~$\Fcal$ is a pair $(Q,\zz)$ 
consisting of a quiver~$Q$ as above together with 
a collection~$\zz$ (an \emph{extended cluster}) 
consisting of algebraically independent (over~$\CC$) 
elements of~$\Fcal$, one for each vertex of~$Q$. 
The elements of~$\zz$ associated with mutable
vertices are called \emph{cluster variables}; they form a
\emph{cluster}. 
The elements associated with the frozen vertices are called
\emph{coefficient variables}. 

For a cluster variable~$z$, 
a~\emph{seed mutation}~$\mu_z$ at the corresponding mutable vertex of~$Q$
transforms $(Q,\zz)$ into the seed $(Q',\zz')=\mu_z(Q,\zz)$ 
defined by $Q'=\mu_z(Q)$ (quiver mutation) and 
$\zz'\!=\!\zz\cup\{z'\}\setminus\{z\}$;  
here~$z'$ is given by the 
\emph{exchange relation}
\begin{equation*}
z\,z'=
\prod_{z\leftarrow y} y
+\prod_{z\rightarrow y
} y\,.
\end{equation*}
(The two products are over the edges directed at and from~$z$, respectively.) 

Seeds related by a sequence of
mutations are called \emph{mutation equivalent}. 
The \emph{cluster algebra $\Acal(Q,\zz)$} associated to a seed
$(Q,\zz)$ is the subring of $\Fcal$
generated by 
all elements of all extended clusters of the seeds mutation-equivalent
to~$(Q,\zz)$. 
%

\section{3. Marked bordered surfaces}

Let $S$ be a connected oriented surface with nonempty
boundary~$\partial S$ and 
finitely many marked points on~$\partial S$,
each of them colored black or white. 
Such a surface is defined by its genus~$g$ together with the patterns of
marked points on its boundary components. 
These patterns are recorded by the \emph{signature}~$\sigma$ 
of the marked surface~$S=S_{g,\sigma}$. 
For example, 
if $g=0$ and $\sigma = [\,\bullet\, | \bullet | \circ \bullet \, \circ\,]$,
then $S_{g,\sigma}$ is a pair of pants 
with two of the boundary components marked by a single black point and
the third component marked by one black and two white points. 

\begin{figure}[ht]
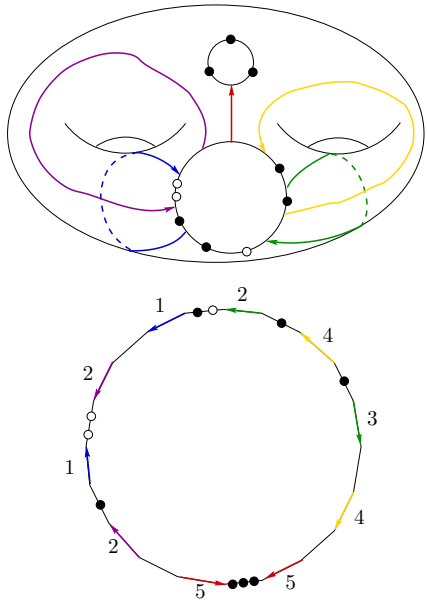

\vspace{-.15in}
\centerline{\scalebox{.8}{\input{sw1a.pstex_t}}}
\vspace{.15in}
\centerline{\scalebox{.8}{\input{sw1b.pstex_t}}}
\vspace{.05in}
\caption{Two-holed surfaces of genus $g=2$, with 
$\sigma=[\,\bullet \bullet \bullet\,|\bullet \bullet \circ \bullet \bullet \circ\,
    \circ\, ]$.
Here $a=3$, $b=7$, $c=5$.
The five cuts produce the 20-gon shown at the bottom.
The labeling indicates which sides are glued to each other.}
\label{fig:sw1a}
\end{figure}

Let us draw several simple non-inter\-secting
curves on~$S$ (``cuts'') such that:
\begin{itemize}
\vspace{-.05in}
 \item $S$ minus the cuts is homeomorphic to a disk;
 \item each cut connects unmarked boundary points; 
 \item for each cut, a choice of direction is made;
 \item each cut is defined up to isotopy that fixes its endpoints.
\end{itemize}
\vspace{-.05in}
Cutting up the surface along these curves yields a polygon
whose sides alternate between boundary and cut segments. 
The perimeter of the polygon represents a walk along boundary components 
and cuts such that the surface always remains on the right side of the walker. 
%
See Figure~\ref{fig:sw1a}. 


If $S$ has $a$ white marked points, $b$~black marked points,
and requires $c$ cuts as above, then we say that $S$ is of
\emph{type}~$(a,b,c)$. 
The type is determined by the genus~$g$ and the
signature~$\sigma$: 
\vspace{-.05in}
\begin{align}
\label{eq:abc-gsigma-1}
&\text{$a$ is the number of white points $\circ$ in $\sigma$;}\\
\label{eq:abc-gsigma-2}
&\text{$b$ is the number of black points $\bullet$ in $\sigma$;}\\
\label{eq:abc-gsigma-3}
&\text{$c=2g-1 +\langle\text{number of boundary components}\rangle$.}
\end{align}

We will connect the combinatorial
topology of a marked surface of type $(a,b,c)$
to the ring of invariants $R_{a,b,c}$
by defining distinguished elements 
in $R_{a,b,c}$ 
associated with particular embedded trivalent~graphs on~$S_{g,\sigma}$
called tensor diagrams.
This is a generalization of the construction 
described in~\cite{fomin-pylyavskyy}; 
in \emph{loc.\ cit.}, the surface $S_{g,\sigma}$ is a disk, so $c=0$ (no cuts). 

\vspace{-.35in}

\begin{figure}[ht]
    \begin{center}
\setlength{\unitlength}{2pt}
\begin{picture}(20,20)(0,0)
\thicklines
\put(10,10){\vector(1,1){10}}
\put(10,10){\vector(-1,1){10}}
\put(10,10){\vector(0,-1){10}}
\end{picture}
\qquad
\begin{picture}(20,20)(0,0)
\thicklines
\put(0,20){\vector(1,-1){10}}
\put(20,20){\vector(-1,-1){10}}
\put(10,0){\vector(0,1){10}}
\end{picture}
\qquad
\begin{picture}(15,20)(0,0)
\thicklines
\put(0,10){\vector(1,0){15}}
\end{picture}
\qquad
\begin{picture}(15,20)(0,0)
\thicklines
\put(0,10){\vector(1,0){15}}
\put(7.5,5){\color{red} \vector(0,1){10}}
\end{picture}
\qquad
    \end{center} 
    \caption{Building blocks for tensor diagrams.
The edges can be bent or stretched.  
}
    \label{fig:webs21}
\end{figure}
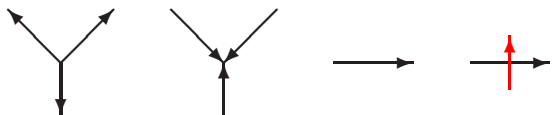 

\section{4. Tensor diagrams on a surface} 

We henceforth assume that $V$ is $3$-dimensional.

First, an informal description. 
Tensor diagrams on a surface $S=S_{g,\sigma}$
are constructed using four
types of building blocks shown in Figure~\ref{fig:webs21}. 
The first three blocks 
represent
three basic $\SL(V)$-invariant tensors:
the volume tensor, the dual volume tensor and the identity tensor. 
The signs of the volumes are well defined
as the orientation of the surface determines a cyclic ordering
on the edges incident to each trivalent vertex. 
The endpoints of each block correspond to the tensor's arguments:
a sink to a vector, a source to a covector. 

The fourth block represents a matrix tensor 
for an element $X \in SL(V)$ associated with the cut shown in red.  
While the matrix tensor is not $\SL(V)$-invariant, it can be used to build
$\SL(V)$-invariant tensors. 
Crossing the red cut in the opposite direction 
represents the matrix tensor for~$X^{-1}$. 

Blocks of these four types can be ``plugged'' into each other (respecting the
orientation) and ``clasped'' at the boundary,
cf.~\cite{fomin-pylyavskyy} for details. 
Since each vertex is a source or a sink, we can replace the
orientation of the edges by a bi-coloring of the vertices
(sinks become black, and sources white).
This leads us to the following definition.

A \emph{tensor diagram} on $S=S_{g,\sigma}$
is a finite bipartite graph $D$ embedded in~$S$, with a fixed 
proper coloring of its vertices into two colors, black and white, 
such that each internal vertex is trivalent,
and each boundary vertex is a marked point of~$S$.
(The embedded edges of~$D$ are allowed to cross each other.)
In addition, $D$ may contain 
oriented loops without vertices. 

We denote by $\operatorname{bd}(D)$ 
(resp., $\operatorname{int}(D)$) the set 
of \emph{boundary} (resp., \emph{internal}) vertices of~$D$. 

Fix a collection of cuts in~$S$, as in Section~3. 
A tensor diagram~$D$ on $S$
defines an $\SL(V)$ invariant
$[D]\in R_{a,b,c}$ obtained by
repeated contraction of elementary $\SL(V)$-invariant tensors. 
The arguments of $[D]$ 
(viewed as a polynomial function on
$(V^*)^a\times V^b \times (\SL(V))^c$)
are interpreted as follows:
\begin{itemize}
\vspace{-.05in}
\item
each black boundary vertex represents a vector argument; 
\item
a white boundary vertex represents a covector argument; 
\item
a cut represents a matrix argument. 
\end{itemize}
\vspace{-.05in}
For the sake of precision, we give a formula for~$[D]$. 
Let $\operatorname{cut}(D)$ denote the set
of points where $D$ crosses the cuts. 
\emph{Edge fragments} are the pieces into which those cuts cut the
edges of~$D$. 
(If an edge is not cut, it forms an edge fragment by itself.) 
Then the invariant $[D]$ is given by
\begin{equation*}
\begin{split}
[D]=&\sum_\ell 
\biggl(\,\prod_{v\in\operatorname{int}(D)}\operatorname{sign}(\ell(v))\biggr)
\biggl(\,\prod_{\substack{v\in\operatorname{bd}(D)\\
\text{$v$ black}}}x(v)^{\ell(v)}\biggr)\\ &\quad
\biggl(\,\prod_{\substack{v\in\operatorname{bd}(D)\\
\text{$v$ white}}}y(v)^{\ell(v)}\biggr)
\biggl(\,\prod_{\substack{v\in \operatorname{cut}(D)\\
}}X_{\ell(v)}\biggr)
\end{split}
\end{equation*}
where 
\begin{itemize}
\vspace{-.05in}
\item
$\ell$ runs over all labelings of the edge fragments in~$D$ by the
numbers~$1,2,3$ such that for each internal vertex~$v$ of~$D$,
the labels of the edges incident to $v$ are 
distinct; 
\item
$\operatorname{sign}(\ell(v))$ is the sign of the 
cyclic permutation defined by the clockwise
reading of those three labels; 
\item
$x(v)^{\ell(v)}$ denotes the monomial $\prod_e x_{\ell(e)}(v)$,
product over all edges $e$ incident to~$v$, and similarly for $y(v)^{\ell(v)}$;
\item 
$X_{\ell(v)}$ denotes the entry $X_{ij}$ of the matrix 
$X \in \SL(V)$ associated with the cut at a vertex $v\in \operatorname{cut}(D)$,
where $i$ and $j$ are the labels of the two adjacent edge fragments.
\end{itemize}
\vspace{-.05in}
A couple of examples are shown in Figure~\ref{fig:sw2}.

\pagebreak[3]

\vbox{

\begin{figure}[ht]
    \begin{center}
\vspace{-.05in}
\scalebox{0.9}{\input{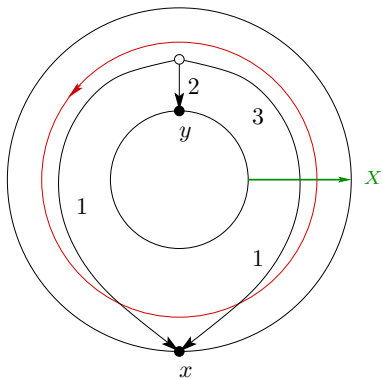}}
    \end{center} 
    \caption{An annulus, with $\sigma=[\bullet|\bullet]$.
A~cut is shown in green. 
The tensor diagram $D_1$ shown in black
defines the invariant $[D_1] = \langle x, X x, y \rangle$. 
The particular proper labeling of edge fragments shown in the picture
contributes the term $- x_1^2 y_2 X_{31}$ to~$[D_1]$.
The red ``floating'' loop $D_2$ defines the invariant $[D_2]=tr(X^{-1})$. 
}
    \label{fig:sw2}
\end{figure}

\vspace{-.15in}

\begin{figure}[hb] 
\scalebox{0.9}{\input{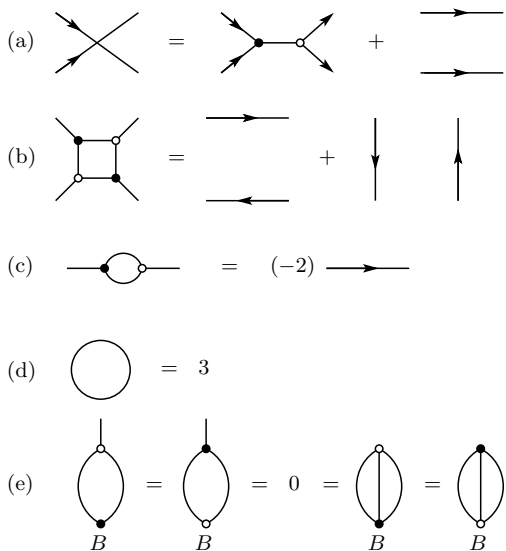}}
\vspace{.05in}
    \caption{Local relations for (invariants associated with)
tensor diagrams. Remember that all
    edges are oriented towards black endpoints. The cycle in
    relation~(d) can be oriented either way.
Marked vertex~$B$ lies on the boundary.}
    \label{fig:skein}
\end{figure}



\begin{figure}[ht]
\vspace{-.2in}
    \begin{tabular}{cl}
\setlength{\unitlength}{1pt}
\begin{picture}(16,20)(0,0)
\put(8,10){\makebox(0,0){(a)}}
\end{picture}
&
\setlength{\unitlength}{1pt}
\begin{picture}(30,20)(0,0)
\thicklines
\put(0,20){\line(3,-2){30}}
\put(0,0){\line(3,2){30}}
\qbezier(0,10)(15,0)(30,10)
\end{picture}
\begin{picture}(16,20)(0,0)
\put(8,10){\makebox(0,0){$=$}}
\end{picture}
\begin{picture}(30,20)(0,0)
\thicklines
\put(0,20){\line(3,-2){30}}
\put(0,0){\line(3,2){30}}
\qbezier(0,10)(15,20)(30,10)
\end{picture}
\\[.2in]
\setlength{\unitlength}{1pt}
\begin{picture}(16,16)(0,0)
\put(8,8){\makebox(0,0){(b)}}
\end{picture}
&
\setlength{\unitlength}{1pt}
\begin{picture}(30,16)(0,0)
\thicklines
\qbezier(0,0)(20,5)(20,10)
\qbezier(30,0)(10,5)(10,10)
\qbezier(10,10)(10,15)(15,15)
\qbezier(20,10)(20,15)(15,15)
\end{picture}
\begin{picture}(16,16)(0,0)
\put(8,8){\makebox(0,0){$=$}}
\end{picture}
\begin{picture}(30,16)(0,0)
\thicklines
\put(0,0){\line(1,0){30}}
\end{picture}
\\[.3in]
\setlength{\unitlength}{1pt}
\begin{picture}(16,10)(0,0)
\put(8,5){\makebox(0,0){(c)}}
\end{picture}
&
\setlength{\unitlength}{1pt}
\begin{picture}(30,10)(0,0)
\thicklines
\qbezier(0,0)(15,15)(30,0)
\qbezier(0,10)(15,-5)(30,10)
\end{picture}
\begin{picture}(16,10)(0,0)
\put(8,5){\makebox(0,0){$=$}}
\end{picture}
\begin{picture}(30,10)(0,0)
\thicklines
\put(0,0){\line(1,0){30}}
\put(0,10){\line(1,0){30}}
\end{picture}
\\[.2in]
\setlength{\unitlength}{1pt}
\begin{picture}(16,20)(0,0)
\put(8,10){\makebox(0,0){(d)}}
\end{picture}
&
\setlength{\unitlength}{1pt}
\begin{picture}(30,20)(0,0)
\thicklines
\put(0,10){\line(1,0){8.8}}
\put(10,10){\circle{2.5}}
\qbezier(11,11)(20,18)(30,0)
\qbezier(11,9)(20,2)(30,20)
\end{picture}
\begin{picture}(16,20)(0,0)
\put(8,10){\makebox(0,0){$=$}}
\end{picture}
\begin{picture}(12,20)(0,0)
\put(5,10){\makebox(0,0){$(-1)$}}
\end{picture}
\begin{picture}(30,20)(0,0)
\thicklines
\put(0,10){\line(1,0){8.8}}
\put(10,10){\circle{2.5}}
\put(30,20){\line(-2,-1){18.7}}
\put(30,0){\line(-2,1){18.7}}
\end{picture}
\\[.2in]
\setlength{\unitlength}{1pt}
\begin{picture}(16,20)(0,0)
\put(8,10){\makebox(0,0){(e)}}
\end{picture}
&
\setlength{\unitlength}{1pt}
\begin{picture}(30,20)(0,0)
\thicklines
\put(0,10){\line(1,0){8.8}}
\put(10,10){\circle*{2.5}}
\qbezier(11,11)(20,18)(30,0)
\qbezier(11,9)(20,2)(30,20)
\end{picture}
\begin{picture}(16,20)(0,0)
\put(8,10){\makebox(0,0){$=$}}
\end{picture}
\begin{picture}(12,20)(0,0)
\put(5,10){\makebox(0,0){$(-1)$}}
\end{picture}
\begin{picture}(30,20)(0,0)
\thicklines
\put(0,10){\line(1,0){8.8}}
\put(10,10){\circle*{2.5}}
\put(30,20){\line(-2,-1){18.7}}
\put(30,0){\line(-2,1){18.7}}
\end{picture}
\end{tabular} 
    \caption{Yang-Baxter-type relations for tensor diagrams.
In relations (a)--(c), edge
    orientations on the left-hand sides can be arbitrary;
the orientations on each right-hand side should match those on the~left.}
    \label{fig:yang-baxter}
\end{figure}
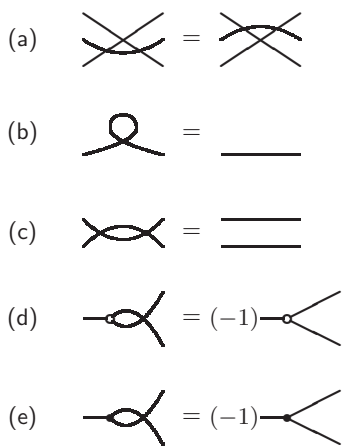
}

\pagebreak[3]

Note that we do allow ``floating'' tensor diagrams that do~not involve any
boundary vertices at all. 
For example, if $D$ is a non-contractible 
oriented loop with~no vertices, 
then $[D]$ is the trace of the 
(non-commutative) product of elements of 
$\SL(V)$ associated with the cuts crossed~by~$D$.  


\begin{theorem}
The invariant $[D]$ defined as above does not change under isotopies
of the cuts and of the tensor diagram~$D$. 
\end{theorem}

The proof consists of checking that local isotopy moves do
not change the invariant. 
We note the importance of the unimodularity condition 
$X\in\SL(V)$ for the matrices associated with the cuts:
in its absence, some local moves 
would contribute a nontrivial factor~$\det(X)$.

\vspace{-.15in}

\section{5. Diagrammatic calculus}
 
The key notational advantage of tensor diagrams is that they naturally
lend themselves to diagrammatic calculations which are more
convenient, less cumbersome, and more intuitive than the more traditional
algebraic formalism. 
This diagrammatic calculus utilizes 
a number of easily verified 
local relations (see Figures~\ref{fig:skein}
and~\ref{fig:yang-baxter})
satisfied by the invariants associated with tensor diagrams.


A  \emph{web} $D$ on an oriented surface $S_{g,\sigma}$ 
is a tensor diagram whose edges do not cross or touch each other, except at
endpoints. 
Each web is considered up to an isotopy of the surface that fixes its
boundary. 
%
The systematic study of webs in case of a disk was pioneered by
G.~Kuperberg~\cite{kuperberg}.

A web is called \emph{irreducible} (or \emph{non-elliptic}) if
it has no contractible loops,  
no pairs of edges enclosing a contractible disk,
and no (unoriented) simple 4-cycles 
whose all vertices are internal 
and which enclose a contractible disk. 

To illustrate, the black tensor diagram in Figure~\ref{fig:sw2} 
is a non-eliptic web, and ditto for the red one
(but not the union of the two, since they cross each other). 
 
An invariant associated to an irreducible web is called a \emph{web
  invariant}. 

One easily checks that by repeatedly applying relations from
Figures \ref{fig:skein}--\ref{fig:yang-baxter}, 
any tensor diagram can be transformed into a (finite) formal
linear combination of non-elliptic webs. 
We call this the \emph{flattening} process. 
The following result
generalizes a theorem by G.~Kuperberg~\cite{kuperberg2}.

\begin{theorem} 
\label{thm:confl}
 The flattening process is confluent:
the formal linear combination of irreducible webs 
that it produces does not depend on the choice of flattening moves.
\end{theorem}

The proof of Theorem~\ref{thm:confl} is based on a diamond lemma argument. 

For a given surface $S_{g,\sigma}$, 
one can define the \emph{skein algebra} $A(S_{g,\sigma})$ of 
(formal linear combinations of) tensor diagrams, 
subject to the relations in Figures  \ref{fig:skein}--\ref{fig:yang-baxter}. 
Multiplication in $A(S_{g,\sigma})$
is given by superposition of diagrams. 
Thus, Theorem~\ref{thm:confl} is a statement about calculations
in~$A(S_{g,\sigma})$. 

\begin{theorem}
\label{thm:surj}
The correspondence $D\mapsto [D]$ extends to a surjective ring homomorphism
$\phi:A(S_{g,\sigma})\to R_{a,b,c}$. 
\end{theorem}

The map $\phi$ is a well-defined ring homomorphism 
since the invariants associated with
tensor diagrams satisfy the relations that define $A(S_{g,\sigma})$. 
Surjectivity follows by checking that each Procesi generator 
of $R_{a,b,c}$ (see Section~1) lies in the image.
Indeed, these generators can be realized by very simple tensor
diagrams, namely bipods, tripods and closed loops. 

\begin{conjecture}
\label{conj:phi-injective}
 The map $\phi$ is injective. 
\end{conjecture}

Conjecture~\ref{conj:phi-injective} can be restated as saying that 
irreducible webs form a linear basis in $A(S_{g,\sigma})$, 
not just a spanning set. 
(In the case of a disk, this was proved by Kuperberg~\cite{kuperberg}.)
This would follow from a verification of the defining relations  
of $R_{a,b,c}$ in the skein algebra
$A(S_{g,\sigma})$. 
This strategy~was successfully implemented in \cite{muller-samuelson}  
in the $\SL_2$ case ($\dim V=2$). 

\pagebreak[3]

\section{6. Special invariants}

In Sections 6--7,
we explain our construction of a
cluster algebra structure in a ring of invariants $R_{a,b,c}$. 
The approach is close to the 
one we used in~\cite[Sections~6--7]{fomin-pylyavskyy} 
in the special case when $S$ is a disk. 
The latter paper provides many examples and pictures which the reader
may find helpful. 

The construction depends on a choice of a marked bordered surface
$S=S_{g,\sigma}$ satisfying
conditions~\eqref{eq:abc-gsigma-1}--\eqref{eq:abc-gsigma-3}. 
We also require that 
\begin{equation}
\label{eq:proper-signature}
\begin{array}{l}
\text{each boundary component~$C$ carries some marked}\\
\text{points, and their colors do not alternate along~$C$.} 
\end{array}
\end{equation}
For a given ring $R_{a,b,c}$, there are typically many choices of 
such a marked surface~$S_{g,\sigma}$.
At least one such choice exists~unless 
$a\!=\!b\!=\!0$, or $a+b\!=\!1$ and $c$ is odd, or
$a\!=\!b\!=\!1$ and $c$ is even.

\begin{figure}[ht]
\vspace{-.1in}
\hspace{-.1in}\scalebox{1}{\input{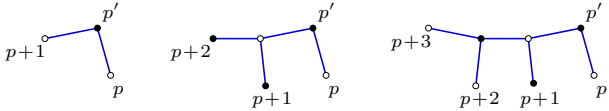}}
\vspace{.1in}
    \caption{Trees $\Lambda_p$. The proxy vertex is denoted by~$p'$.}
    \label{fig:caterpillars}
\end{figure}

Along each boundary component, say with $m$ marked points,
we number those marked points by the elements of~$\ZZ/m\ZZ$. 
For each component, we are free to decide between the clockwise 
and counterclockwise order;
the construction works for any choices.
For simplicity, let us assume that each component is numbered 
counterclockwise.
(More precisely, the labels increase as we move along the boundary so
that the surface remains on the right.) 
For each boundary vertex~$p$,  
we define two trees $\Lambda_p$ and $\Lambda^p$ embedded in~$S$. 
If~$p$ is black, then $\Lambda_p$ has one vertex,
namely~$p$, and no edges.
If~$p$ is white, then place a new black vertex 
(the \emph{proxy} of~$p$) next to~$p$ inside $S$, and connect it to~$p$. 
Then look at $p+1$. If $p+1$ is white,
then connect it to the proxy vertex, and stop; 
see Figure~\ref{fig:caterpillars} on the left.
In general, we proceed clockwise from~$p$ until we find two
consecutive vertices of the same color (here we need
the non-alternating condition~\eqref{eq:proper-signature}), 
then build a caterpillar-like
bi-colored tree as shown in Figure~\ref{fig:caterpillars}.
The graph $\Lambda^p$ (and the proxy vertex in the case when $p$ is white) 
is defined in the same way, with the colors swapped. 

\vspace{-.2in}

\emph{Arcs} on~$S$ are simple curves connecting marked points,
considered up to an isotopy fixing the endpoints of a curve. 
Contractible loops are not allowed.
In contrast with~\cite{cats1},
we allow curves isotopic to a boundary segment between consecutive
marked points.
We moreover identify such an arc with the corresponding boundary segment. 

We next define three families of invariants
associated with configurations of arcs: 

\emph{Special invariants $J_p^q(\alpha)$.} 
Let $\alpha$ be an arc with endpoints $p$ and~$q$. 
Then $J_p^q(\alpha)$ is the invariant 
defined by the tensor diagram obtained by connecting the trees
$\Lambda_p$ and~$\Lambda^q$ by a single edge $e$ running along the arc~$\alpha$.
At one end, $e$~connects to $p$ if the latter is black,
otherwise to the proxy of~$p$.
At the other end, $e$~connects to $q$ if the latter is white, 
otherwise to its proxy. 
Whenever $e$ terminates at a proxy vertex, it approaches the latter from
the direction opposite from the nearest boundary. 

\emph{Special invariants $J_{pqr}(\alpha \beta \gamma)$ and
$J^{pqr}(\alpha \beta \gamma)$.} 
Let $p,q,r$ be marked points connected pairwise by the 
arcs $\alpha,\beta,\gamma$ cutting out a
(contractible) \emph{triangle} $\alpha \beta \gamma$.
The invariant $J_{pqr}(\alpha \beta \gamma)$
is defined by the tensor diagram obtained as follows. 
Place a white vertex $s$ in the middle
of the triangle $\alpha \beta \gamma$, and 
connect $s$ inside the triangle to each of $\Lambda_p, \Lambda_q,
\Lambda_r$. 
As the other endpoints of these three edges,
use the vertices $p,q,r$ whenever they are black, otherwise take
their respective proxies. 
For $J^{pqr}(\alpha \beta \gamma)$, reverse the roles of the colors. 
%
To give an example,
the black tripod web in Figure~\ref{fig:sw2} is the special invariant
$J_{xyx}(\alpha \beta \delta)$, 
in the notation of Figure~\ref{fig:sw34}.

\emph{Special invariants $J_{pq}^{rs}(\alpha \beta)$.}
Let $p,q,r,s$ be four marked points. 
Let $\alpha$ be an arc connecting $p$ to~$r$, 
and let $\beta$ be an arc connecting $q$ to~$s$, 
so that $\alpha$ and $\beta$ cross exactly once.
The invariant $J_{pq}^{rs}(\alpha \beta)$ is defined by the tensor diagram
obtained as follows. 
Place a white vertex~$W$ and a black vertex~$B$ near the 
intersection of $\alpha$ and~$\beta$. 
Connect $W$ and $B$ by an edge. 
Connect $W$ to $L_p$ and~$L_q$ along $\alpha$ and~$\beta$,
respectively, 
using $p$ and $q$ if these two are
black, or else using their (black) proxies as needed. 
Similarly, connect $B$ to the appropriate white vertices 
in $L^r$ and~$L^s$ along $\alpha$ and~$\beta$. 
Make sure that the five edges incident to $B$ and $W$ do not cross each other. 



The following four propositions generalize their respective counterparts
in~\cite[Section~6]{fomin-pylyavskyy}. 

\begin{proposition}
\label{prop:special=0}
If $p$ and~$q$ are not adjacent on the same boundary component, 
then $J_p^q(\alpha)$ is a nontrivial invariant. 
If $p$ and~$q$ are adjacent, and $\alpha$ is isotopic to 
a boundary segment between $p$ and~$q$,
then exactly one of the two invariants
$J_p^q(\alpha)$ and $J_q^p(\alpha)$ vanishes.
Specifically, if $p$ is white, then $J_p^{p+1}(\alpha)=0$;
if $p$ is black, then $J_{p+1}^p(\alpha)=0$. 
\end{proposition}

\begin{proposition}
\label{prop:specweb}
Any nonzero special invariant is a web invariant
(i.e., an invariant defined by an irreducible web). 
\end{proposition}

We call a nonzero special invariant \emph{indecomposable} if it
does not factor as a product of two or more special invariants. 

\begin{proposition}
\label{prop:special-inv-factoring}
Any nonzero special invariant is represented uniquely as a product of
 indecomposable special invariants. 
\end{proposition}

There is a simple algorithm for finding such a factorization called 
{\it {arborization}}, cf.\ Conjecture \ref{conj:arb}.

Special invariants are \emph{compatible} if their product is a
single web invariant. 
A special invariant compatible with any other special invariant is
called a \emph{coefficient invariant}.
This terminology anticipates the
appropriate notions of (a) compatibility of cluster variables
and (b) coefficient variables. 


\begin{proposition}
\label{prop:coeff-special}
On a surface $S=S_{g,\sigma}$ of type $(a,b,c)$, 
there are $a+b$ coefficient invariants---unless $S$ is a disk with
$a+b\le 4$. 
These $a+b$ invariants are 
the (indecomposable) nonzero invariants of the form $J_p^{p \pm 1}$.
Moreover the product of any of them and any web invariant is
 a web invariant. 
\end{proposition}


To illustrate, the yellow tensor diagram in Figure~\ref{fig:sw34} 
is one of the two coefficient special invariants,
namely $J_{y}^y(\gamma)$. The other coefficient is $J_{x}^x(\alpha)$.

To obtain exchange relations for our cluster
algebras, we will need certain 
$3$-term skein relations for special invariants. 

\begin{proposition}
\label{prop:3-term-1}
Let $\alpha\beta\gamma$ be a triangle formed by arcs 
$\alpha, \beta, \gamma$ which respectively connect $(p,q)$,
$(q,r)$, and $(p,r)$. 
Then
\begin{align}
\notag
J_{pqr}(\alpha \beta \gamma) J^{pqr}(\alpha \beta \gamma) =
&J_r^p(\gamma) J_q^r(\beta) J_p^q(\alpha)\\ 
\label{eq:3-term-1}
+ &J_r^q(\beta) J_p^r(\gamma) J_q^p(\alpha)\,.
\end{align}
\end{proposition}

\begin{figure}[ht]
    \begin{center}
\input{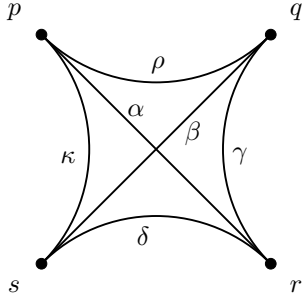}
\end{center}
\caption{Arcs $\gamma$, $\delta$, $\kappa$, and~$\rho$ stay close to
  the arcs $\alpha$ and~$\beta$, 
ensuring that the quadrilateral $\gamma\delta\kappa\rho$ is
  contractible.}
    \label{fig:sw14}
\end{figure}

\begin{proposition}
Let $p,q,r,s$ be marked points.
Let $\alpha, \beta, \gamma, \delta, \kappa, \rho$
be arcs connecting them as in Figure~\ref{fig:sw14}.
Then
\begin{align}
\label{eq:3-term-2}
J_p^r(\alpha) J_{qrs}(\beta \gamma \delta) &= J_q^r(\gamma)
J_{prs}(\alpha \delta \kappa)+ J_s^r(\delta) J_{pqr}(\alpha \rho \gamma)\,;\\
\label{eq:3-term-3}
J_p^r(\alpha) J_s^q(\beta)
&= J_s^r(\delta) J_p^q(\rho) + J_{sp}^{qr}(\alpha \beta)\,;\\
\notag
J_p^r(\alpha) J_{rs}^{pq}(\alpha \beta) &=J_p^q(\rho)
J_r^p(\alpha) J_s^r(\delta)\\
\label{eq:3-term-4}
&+ J_{prs}(\alpha \delta \kappa) J^{pqr}(\alpha \rho \gamma)\,;\\
\notag
J_r^p(\alpha) J_{sp}^{qr}(\alpha \beta) &= J_r^q(\gamma)
J_s^p(\kappa) J_p^r(\alpha)\\
\label{eq:3-term-5}
&+ J_{prs}(\alpha \delta \kappa)J^{pqr}(\alpha \rho \gamma)\,.
\end{align}
\end{proposition}

In each of the identities \eqref{eq:3-term-1}--\eqref{eq:3-term-5}, 
some special invariants might vanish; others might factor further.
Once everything is expressed in terms of indecomposable special invariants, 
one either gets a tautological formula $A=A$, 
or else a genuine \hbox{$3$-term} relation. 
We call the latter relation the \emph{distilled form} of the original
one. Examples can be found in \cite{fomin-pylyavskyy}.


\section{7. Seeds associated with triangulations}

We next describe how to associate a seed to a triangulation~$T$ 
of a surface $S=S_{g,\sigma}$ by arcs connecting marked points
on~$\partial S$. 

We begin by building a collection $K(T)$ of special invariants, as
follows:  
\begin{itemize}
\vspace{-.05in}
\item 
for each arc $\alpha$ in $T$
connecting some pair of marked points $p$ and $q$, 
include $J_p^q(\alpha)$ and $J_q^p(\alpha)$ in $K(T)$;
\item 
for each triangle $\alpha\beta\gamma$ in~$T$ 
with vertices $p,q,r$ (listed clockwise), 
include $J_{pqr}(\alpha \beta \gamma)$.
\end{itemize}
\vspace{-.05in}
The \emph{extended cluster} $\zz(T)$ associated with~$T$ 
is the set of all indecomposable special invariants which appear in
factorizations of nonzero elements of $K(T)$ into indecomposables. 
The \emph{cluster} $\xx(T)$ 
consists of all non-coefficient invariants in~$\zz(T)$. 
An~example is shown in Figure \ref{fig:sw34}.

\begin{theorem}
\label{th:z(T)}
For any triangulation~$T$ as above, 
\vspace{-.05in}
\begin{itemize}
\item
$\zz(T)$ contains the entire set of
coefficient invariants; 
\item
$\zz(T)$ consists of $3(a+b)+8(c-1)$ invariants;
\item
all special invariants in $\zz(T)$ are pairwise compatible. 
\end{itemize}
\end{theorem}


Our next goal is to write the exchange relations for an \linebreak[3]
extended cluster~$\zz(T)$. 
Encoding these relations by an appropriate quiver will then 
complete the construction of a seed  associated with a triangulation~$T$. 

\begin{proposition}
\label{pr:thin-quads}
Let $(p,p+1,p+2,s)$ be four distinct marked points, the first three of
them consecutive on the same boundary component.
Let $\alpha, \beta, \gamma, \delta, \kappa, \rho$ be the arcs of $T$
connecting the four points as in Figure~\ref{fig:sw14},
with $q = p+1$, $r = p+2$.
\begin{align}
\notag
&\text{If $p$ is white and $p+1$ is black,
then} \\
\label{eq:3-term-3wb}
&\qquad J_{p+2}^p(\alpha) J_{p+1}^s(\beta) = J_{p+1}^p(\rho)
J_{p+2}^s(\delta) + J_p^s(\kappa).
\end{align}
\vspace{-.05in}
\begin{align}
\notag
&\text{If $p$ is black and $p+1$ is white,
then}\\
\label{eq:3-term-3bw}
&\qquad
J_p^{p+2}(\alpha) J_s^{p+1}(\beta) = J_p^{p+1}(\rho) J_s^{p+2}(\delta)
+ J_s^p(\kappa).
\end{align}
\end{proposition}

A side of a triangle 
is called \emph{exposed} if it lies on the
boundary of~$S$, 
connecting two adjacent marked points. 

We are now prepared to generate the exchange relations for~$\zz(T)$. 
Let us write the following identities:
\vspace{-.03in}
\begin{itemize}
\item
for each triangle $\alpha\beta\gamma$ of~$T$, 
write formula~\eqref{eq:3-term-1};
\item
for each diagonal $pr$ in~$T$ separating triangles $pqr$ and $prs$:
\begin{itemize}
\vspace{-.03in}
\item 
write formula~\eqref{eq:3-term-2};
\item 
if one of the sides of $pqr$ is exposed, 
write~\eqref{eq:3-term-4}--\eqref{eq:3-term-5};
\item 
if two sides of $pqr$ are exposed, 
write the appropriate instance of~\eqref{eq:3-term-3wb}
or~\eqref{eq:3-term-3bw}, if applicable. 
\end{itemize}
\end{itemize}
One can check
that 
each of the two monomials on the right-hand side of each of resulting
relations will either vanish (in which case we discard the relation) 
or else factor into nontrivial indecomposables. 
In the latter case, we distill the relation
to obtain a $3$-term relation involving
indecomposable special invariants all of which, with the exception on
the second factor on the left, belong~to~$\zz(T)$. 
Some of the relations obtained by this procedure 
may be identical to each other. 

To obtain the final list of exchange relations,
we should also inspect all instances where 
there is another triangulation~$T'$ with
the same cluster~$\xx(T')=\xx(T)$
(as in \cite[Proposition~7.4]{fomin-pylyavskyy}),
then check whether applying the above recipe to~$T'$ yields any
additional $3$-term relations,
cf.\ \cite[Definition~7.8]{fomin-pylyavskyy}. 

\begin{proposition}
\label{prop:special-seeds-exchanges}
The procedure described
above yields 
one relation of the form 
$x x'=M_1+M_2$ for each $x\in\xx(T)$;
here $M_1,M_2$ are monomials in the elements of~$\zz(T)$. 
\end{proposition}

The \emph{quiver $Q(T)$ associated
  with a triangulation~$T$} is defined by the relations
obtained via the procedure outlined above. 
An example is shown in Figure~\ref{fig:sw34} on the right. 
To be precise, the relations define $Q(T)$ up to simultaneous reversal of
  direction of all edges incident to any subset of connected components
  of the mutable part of the quiver. 
(Usually there will be a~single connected component.)
This ambiguity in the definition of $Q(T)$ does not create any
problems since the choices involved do not affect 
the actual cluster structure. 

\begin{figure}[h]
\begin{center}
\scalebox{.8}{
\input{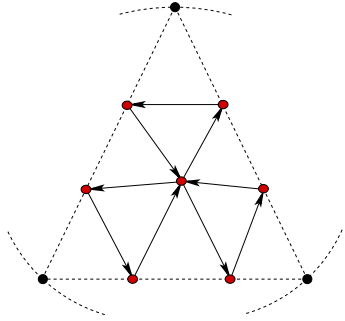}
}
\end{center}
\caption{A building block for a triangle without exposed sides.}
    \label{fig:sw5a}
\end{figure}

\begin{figure}[h]
\vspace{-.1in}
\begin{center}
\input{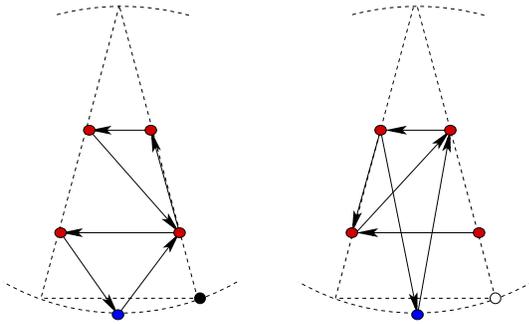}
\vspace{.1in}
\end{center}
\caption{Building blocks for triangles with one exposed side.}
    \label{fig:sw5b}
\end{figure}

In many cases, the quiver $Q(T)$
for a triangulation cluster in $\Acal(S_{g,\sigma})$ 
can be assembled by gluing together 
the building blocks associated to individual triangles in~$T$.
Most common building blocks are shown in
Figures~\ref{fig:sw5a}--\ref{fig:sw5b}. 

\begin{figure*}[b]
\vspace{-.1in}
    \begin{center}
\hspace{-.5in}
\scalebox{0.9}{\input{sw3a.pstex_t}}
\qquad\qquad \scalebox{0.9}{\input{sw34a.pstex_t}}
    \end{center} 
\caption{An annulus $S=S_{g,\sigma}$ with $g=0$ and
  $\sigma=[\bullet|\bullet]$, cf.\ Figure~\ref{fig:sw2}. 
Arcs $\alpha$ and $\gamma$ are boundary segments;
arcs $\beta$ and $\sigma$ cut~$S$ into triangles $\alpha\beta\sigma$
and $\beta\sigma\gamma$.\\
Shown in the middle is the cluster obtained from this triangulation.
Shown on the right is the quiver~$Q(T)$.
This is the quiver of a $Q$-system of type~$A_2$,
cf.\ e.g.\ \cite{difrancesco-survey}.\\
Mutation in the $Q$-system direction corresponds to the action of the 
Dehn twist. 
The red floating loop from Figure~\ref{fig:sw2}, when expressed in terms
of this cluster, is given by\\[.05in]
\centerline{$\displaystyle[D] = \dfrac{J_y^x(\beta)^2 J_x^y(\beta) J_x^y(\delta) +
  J_y^x(\beta) J_y^x(\delta) J_x^y(\beta)^2+ 
 J_y^x(\beta)^2 J_x^x(\alpha) J_y^x(\delta) + J_x^y(\beta)
 J_y^y(\gamma) J_x^y(\delta)^2+ J_y^x(\beta)^2 J_x^y(\delta)^2} 
{J_y^x(\beta) J_x^y(\beta) J_y^x(\delta) J_x^y(\delta)}$.}
} 
    \label{fig:sw34}
\end{figure*}

\pagebreak[3]

\vspace{-.2in}


\section{8. Main results}

We now state our main result: the construction presented above does
indeed yield a desired cluster structure in~$R_{a,b,c}\,$. 


\begin{theorem}
\label{th:main}
For any triangulation~$T$ of the surface~$S_{g,\sigma}$,
the pair~$(\zz(T),Q(T))$ 
constructed as above forms a seed in the field of fractions for
the ring of invariants~$R_{a,b,c}$. 
The corresponding cluster algebra $\Acal(S_{g,\sigma})$ lies
inside~$R_{a,b,c}$. 
This cluster structure 
does not depend on the choice of triangulation~$T$.
\end{theorem} 

The proof of Theorem~\ref{th:main} follows the same strategy 
as that of 
\cite[Theorem~8.1]{fomin-pylyavskyy}.
The key steps of the proof are: 
\begin{itemize}
\vspace{-.05in}
\item
a combinatorial verification that seeds
$(\zz(T),Q(T))$ associated to different triangulations are mutation
equivalent; 
\item
an argument based on the algebraic
Hartogs' principle \cite[Corollary~3.7]{fomin-pylyavskyy},
which uses the fact that the rings $R_{a,b,c}$ are unique
factorization domains, see \cite[Theorem~3.17]{popov-vinberg}; and
\item
showing that indecomposable special invariants are irreducible.
\end{itemize}
\vspace{-.05in}
In the case when $S$ is a disk (equivalently, $c=0$) treated
in~\cite{fomin-pylyavskyy}, we proved a stronger claim 
$\Acal(S_{g,\sigma})=R_{a,b,c}$
using the fact that all Weyl generators of
$R_{a,b,c}$ are special invariants---hence lie in
$\Acal(S_{g,\sigma})$. 
This argument does not work for $c>0$, as the special
invariants do not generate the ring~$R_{a,b,c}$.

Let $\Ucal(S_{g,\sigma})$ denote the \emph{upper cluster algebra}~\cite{ca3} 
associated with~$\Acal(S_{g,\sigma})$. 

\begin{theorem} 
\label{thm:upper}
 Any web invariant, when expressed in terms of any seed,
is given by a Laurent polynomial. 
Thus
\begin{equation}
\label{eq:sandwich}
\Acal(S_{g,\sigma}) \subseteq R_{a,b,c} \subseteq
\Ucal(S_{g,\sigma}).
\end{equation}
\end{theorem}

\begin{figure}[h]
\vspace{-.2in}
    \begin{center}
\input{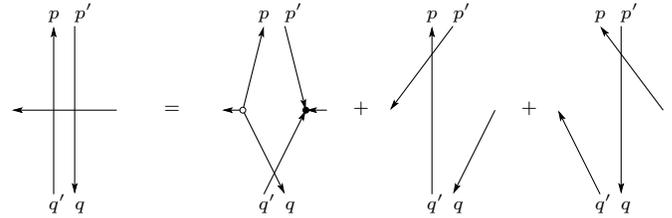}
\vspace{.05in}
\end{center}
\caption{Local relation used in the proof of Theorem~\ref{thm:upper}.}
    \label{fig:sw13}
\end{figure}

We sketch a proof of Theorem~\ref{thm:upper}. 
In view of Theorem~\ref{thm:surj}, 
it suffices to show that any web invariant $[D]$ lies
in $\Ucal(S_{g,\sigma})$.
According to \cite[Theorem 1.5]{ca3} it suffices to check the Laurent
condition with respect to some seed
together with all the seeds obtained from it by a
single mutation. 
We shall explain how to establish Laurentness with respect 
to a seed associated with a triangulation~$T$.
(The verification for the adjacent seeds can be done in a similar fashion.)
We need to show that by repeatedly multiplying $[D]$ by 
elements of $\zz(T)$, we can obtain 
a linear combination of elements of~$\zz(T)$. Let
$\alpha$ be an arc in $T$ with endpoints $p$ and~$q$. The idea is to
multiply $[D]$ by $J_p^q(\alpha) J_q^p(\alpha)$
sufficiently many times for the result to become compatible with
both $J_p^q(\alpha)$ and~$J_q^p(\alpha)$. 
To achieve that, we repeatedly use
the local relation in Figure~\ref{fig:sw13} to get rid of
the crossings between $D$ and~$\alpha$. Once this is done, all the
resulting webs are going to be contained inside individual triangles
of~$T$. This can be shown to imply that they factor into
special invariants in $\zz(T)$ and/or invariants of the form
$J^{pqr}(\alpha\beta\gamma)$, for $\alpha\beta\gamma$ a triangle
in~$T$.
The claim then follows from Proposition~\ref{prop:3-term-1}. 

\vspace{-.1in}

\begin{conjecture}
Under localization of coefficient variables,
all three rings appearing in \eqref{eq:sandwich} become equal. 
\end{conjecture}

When the coefficients are not localized, there is generally
a gap between the cluster algebra $\Acal(S_{g,\sigma})$ and the upper
cluster algebra $\Ucal(S_{g,\sigma})$. 
This is already true in the $\SL_2$ case, and 
can be checked using the model in~\cite{cats1}. 

\pagebreak[3]

Cluster algebras $\Acal(S_{g,\sigma})$ behave in a functorial way
under two types of embeddings, cf.\ \cite[Theorems
  8.8--8.9]{fomin-pylyavskyy}. 
The precise statements 
(see Theorems~\ref{th:drop-vertex}--\ref{th:fork} below) 
make use of the notion of \emph{cluster
  subalgebra}, see \cite[Definition~8.7]{fomin-pylyavskyy}. 

\begin{theorem}
\label{th:drop-vertex}
Let $\sigma$ and $\sigma'$ be two signatures for the same
surface with boundary such that 
{\rm (i)} both $\sigma$ and $\sigma'$ satisfy \eqref{eq:proper-signature} 
and {\rm (ii)} $\sigma'$ is obtained from $\sigma$ by removing a single
marked point. 
Let $R_{a,b,c}$ and $R_{a',b',c}$ 
be the corresponding rings of invariants. 
Then the image 
of $\Acal(S_{g,\sigma'})$ under the natural
embedding $R_{a',b',c}\to R_{a,b,c}$
is a cluster subalgebra of~$\Acal(S_{g,\sigma})$. 
\end{theorem}


\begin{theorem}
\label{th:fork}
Let $\sigma$ and $\sigma'$ be two signatures for the same
surface with boundary such that 
{\rm (i)} both $\sigma$ and $\sigma'$ satisfy \eqref{eq:proper-signature} 
and {\rm (ii)} $\sigma'$ is obtained from $\sigma$ by 
replacing two consecutive entries of the same color by a single entry
of the opposite color. 
Interpreting this operation algebraically as a cross product, 
consider the corresponding embedding
$R_{\sigma'}\to R_\sigma$. 
The image of $\Acal(S_{\sigma'})$ under this embedding 
is a cluster subalgebra of~$\Acal(S_{g,\sigma})$. 
\end{theorem}

The proofs of Theorems \ref{th:drop-vertex}--\ref{th:fork}
are similar to those of \cite[Theorems 8.8--8.9]{fomin-pylyavskyy}.

\begin{corollary}
\label{cor:planar-tree}
Let $D$ be a web in $S_{g,\sigma}\,$. 
If $D$ is a tree (without clasped leaves)
whose leaves lie on at most three boundary components,
then the web invariant $[D]$ is a cluster or coefficient
variable in~$\Acal(S_{g,\sigma})$.
\end{corollary}

\vspace{-.1in}

\section{9. Examples}

\begin{figure}[ht]
\vspace{-.2in}
\begin{center}
\input{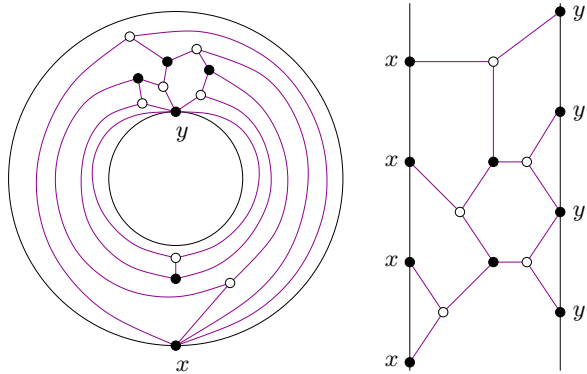}
\vspace{.1in}
\end{center}
\caption{
Two representations of the same cluster variable by tensor diagrams.\\
Left: a non-elliptic web. Right: a tree drawn on the universal cover.}
    \label{fig:sw6}
\end{figure}

Let us review the example in Figure~\ref{fig:sw34},
cf.\ also Figures~\ref{fig:sw2} and~\ref{fig:sw6}. 
Here the surface $S_{g,\sigma}$ is an annulus with one
black marked point on each boundary component.
Thus $g=0$, $\sigma=[\bullet|\bullet]$, $a=0$, $b=2$, $c=1$. 
Looking at the quiver $Q(T)$ for the triangulation of $S_{g,\sigma}$ 
shown in Figure~\ref{fig:sw34},
we recognize that the \emph{cluster type} of
$\Acal(S_{g,\sigma})$ 
(i.e., the mutation equivalence class of the mutable part of~$Q(T)$)
is that of the \emph{$Q$-system} of type~$A_2$. 
For background on $Q$-systems and their connections to cluster
algebras, see, e.g., \cite{difrancesco-survey} and references therein.

\pagebreak[3]

\begin{figure}[h!]
\vspace{-.2in}
    \begin{center}
\input{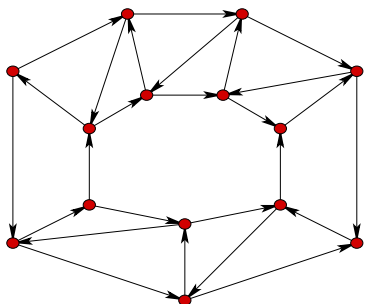}
\end{center}
\caption{The cluster type of an annulus 
with 
$\sigma = [\bullet \bullet \bullet \bullet | \bullet \bullet \,
  \bullet]$.}
    \label{fig:sw9}
\end{figure}

\vspace{-.15in}
\begin{figure}[h!]
    \begin{center}
\input{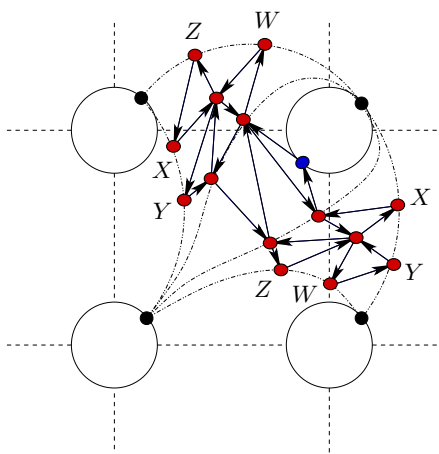}
\end{center}
\caption{A triangulation~$T$ of a one-holed torus 
with a single black marked point,
drawn inside a fundamental domain of the universal cover. 
After opposite sides of the domain are glued, 
some pairs of vertices of the quiver get identified as shown.
Thus $Q(T)$ has 10 mutable vertices.}
    \label{fig:sw8}
\end{figure}


Figure~\ref{fig:sw6} shows a cluster variable in the same cluster
algebra obtained from the seed in
Figure~\ref{fig:sw34} by mutating 
away from the $Q$-system direction. 

More generally, we expect that the cluster type of $Q$-systems of
type~$A_n$ 
arises from a natural cluster structure on the ring of invariants
$\CC[V^2 \times \SL(V)]^{\SL(V)}$
where $\dim V=n+1$.

\smallskip

We next describe a general recipe for determining the cluster type of
$\Acal(S_{g,\sigma})$ 
in the case when $S_{g,\sigma}$ is an annulus 
whose boundary components carry $b_1$ and $b_2$ marked
points, respectively, all of them black. 
(Thus $R_{a,b,c}=R_{a,b_1+b_2,1}$ is the ring of $\SL(V)$-invariants of
$b_1 + b_2$ vectors in~$V\cong\CC^3$ and one matrix in~$\SL(V)$.) 
Take $b_1 + b_2$ squares and glue them together into an annular strip. 
Draw ``parallel'' diagonals inside $b_1-1$ consecutive squares, 
then skip a square, 
then draw parallel diagonals in $b_2-1$ squares. 
Orient the edges of this graph
so that all triangles and unfilled quadrilaterals become oriented
$3$- and $4$-cycles. 
The resulting oriented graph is the mutable part of~$Q(T)$,
for some triangulation of the annulus~$S_{g,\sigma}$. 
An example with $b_1 = 4$ and $b_2 = 3$ is shown in Figure~\ref{fig:sw9}. 

\smallskip

Figure~\ref{fig:sw8} shows the simplest example 
in which $S_{g,\sigma}$ is a surface of positive genus.


\section{10. Conjectures}

The \emph{unclasping} of a tensor diagram~$D$ is the graph obtained
from~$D$ by splitting every boundary vertex~$p$, say of
degree~$k$, into $k$ distinct vertices serving as 
endpoints of the edges formerly incident to~$p$. 
This operation corresponds to \emph{polarization} of invariants. 
We call a tensor diagram~$D$
whose unclasping has no cycles a \emph{forest diagram}; 
if moreover the unclasping is connected, we call $D$ a \emph{tree diagram}. 
Please note that such a diagram~$D$ does not have to be planar. 

\begin{conjecture}[{\rm cf.\ \cite[Conjectures 9.3, 10.1]{fomin-pylyavskyy}}] 
\label{conj:web-and-tree}
All cluster 
variables are web invariants. 
A web invariant $z$ is a cluster monomial 
if and only if $z = [D]$
for some forest diagram~$D$.
\end{conjecture}

The cluster variable shown in Figure~\ref{fig:sw6} is a good 
supporting example for Conjecture~\ref{conj:web-and-tree}:
it can be represented either by a web or by a tree.
Note that the tree form is self-intersecting on the
surface but not on the universal cover. 
Generally speaking, Conjecture~\ref{conj:web-and-tree}
allows for the possibility that the tree form has crossings 
even after being lifted to the universal cover.

\pagebreak[3]

\begin{conjecture}[{\rm cf.\ \cite[Conjecture 10.6]{fomin-pylyavskyy}}] 
\label{conj:arb}
The forest form of a cluster monomial can be found from the web form
 via the \linebreak[3]
\emph{arborization} 
algorithm described in~{\rm \cite[Section~10]{fomin-pylyavskyy}}.
\end{conjecture}

\noindent
\textbf{\small Exercise.}
Check that arborization of the web shown~in 
Figure~\ref{fig:sw6} on the left 
yields the tree tensor diagram on the right.

\begin{conjecture}[{\rm cf.\ \cite[Conjecture 9.2]{fomin-pylyavskyy}}]
\label{conj:compat}
Two cluster (or coefficient) variables are compatible if and only if
their product is a web invariant. 
\end{conjecture}

\noindent
\textbf{\small Exercise.}
Verify that the product of any two of the four cluster 
variables in Figure~\ref{fig:sw34} can be represented by a single web.

\begin{conjecture}[{\rm cf.\ \cite[Conjecture 8.12]{fomin-pylyavskyy}}]
 Reversal of direction 
in the definition of special invariants 
on any subset of boundary components 
does not change the cluster structure $\Acal(S_{g,\sigma})$.
\end{conjecture}

\begin{conjecture}
 Changing all colors of marked points on any subset of boundary
 components does not change the cluster type of $\Acal(S_{g,\sigma})$.
\end{conjecture}

\noindent
\textbf{\small Exercise.}
Verify that the cluster type of
 $\Acal(S_{g,\sigma})$ for 
$\sigma = [\bullet\, | \circ\,]$ is the 
same as in the case $\sigma = [\bullet\, | \bullet\,]$
(the $Q$-system of type~$A_2$). 

\vspace{-.1in}

\section{11. Remarks on additive bases}

A \emph{cluster monomial} is a monomial in the elements of
any extended cluster. 
Cluster monomials are expected to appear in all
``important'' additive bases of cluster algebras 
and closely related rings (skein algebras, 
upper cluster algebras). A~key challenge is to describe the elements
of those bases which are not cluster monomials. 
In the case of cluster algebras associated with surfaces~\cite{cats1}, 
the role of such additional elements is essentially played by floating
$SL_2$ webs. 
We expect a similar situation in the setting of this paper: 
the extra elements should come from the webs which
do not arborize to a forest tensor diagram.

Comparisons of Kuperberg's basis of web invariants 
with Lusztig's \emph{canonical basis} were made
in \cite{khovanov-kuperberg, robert}. 
See \cite{leclerc-icm, musiker-schiffler-williams} for additional
details and references.
As noted in
\cite{fomin-pylyavskyy}, the web basis might coincide
with Lusztig's \emph{semicanonical basis} whenever the
latter is defined.


Our approach meshes well with the philosophy of \cite{BZ, BZ2} 
according to which two elements of 
a ``canonical'' additive basis are compatible if and only if they
quasi-commute.
There is a quantized version of $A_2$
spider relations~\cite{kuperberg2} 
which distinguishes between two ways
of making two strands of a tensor diagram cross:
one chooses which strand goes above the other. 
Whenever the product of webs is a single web (suggesting compatibility,
cf.\ Conjecture~\ref{conj:compat}), 
different crossing patterns yield elements of the quantized skein
algebra that
differ by a scalar factor of the form~$q^m$---thus the original webs
quasi-commute. 

\vspace{-.1in}

\section{12. Fock-Goncharov cluster algebras
}

In their groundbreaking work on higher Techm\"uller theory,
V.~Fock and A.~Goncharov \cite{fg-ihes} 
introduced a family of cluster algebras
which depend on a marked surface~$S$ 
and a semisimple Lie group~$G$. 
For $G=\SL(V)$, 
their construction produces cluster structures in the rings 
of $\SL(V)$ invariants of collections 
of elements of $\SL(V)$
and \emph{affine flags} in~$V$. 
We next describe a modifi\-cation of the tensor diagram
calculus, and of the main construction of this paper, 
that naturally gives the Fock-Goncharov
cluster algebras in the case $G=\SL_3$. 

Let $\dim V=3$. 
Then an affine flag in $V$ is nothing but a vector-covector pair $(v,w)$ with
$\langle v,w\rangle=0$, $v\neq 0$, $w\neq 0$. 
Let us plant such a pair $(v,w)$ 
at each marked point on~$\partial S$. 
Combinatorially, this is encoded by letting every marked point 
carry both black and white colors, 
making it \emph{bivariant} (i.e., either covariant or contravariant as
we please). 
We then define tensor diagrams and associated invariants 
just as before, except that instead of proxies, 
we simply take either the vector or the covector from the corresponding
affine flag, as needed. 
 
Given a triangulation of~$S$, we build an extended cluster
by including all invariants $J_p^q(\alpha)$ 
and $J_{pqr}(\alpha \beta \gamma)$ as in Section~7. 
The corresponding quiver is then constructed using the rule in
Figure~\ref{fig:sw5a}. 
There is no need to worry about
factoring our invariants into indecomposables: 
just use special invariants as cluster or coefficient variables. 

\begin{theorem}
 The construction described above reproduces the triangulation 
seeds of the Fock-Goncharov cluster algebra~{\rm \cite{fg-ihes}}. 
\end{theorem}

\begin{figure}[h]
\vspace{-.2in}
    \begin{center}
\input{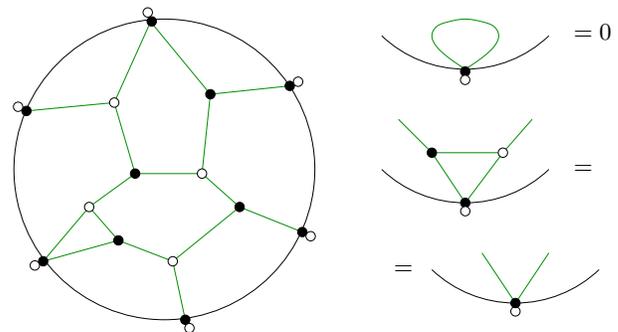}
\end{center}
\caption{Left: a bi-web. This bi-web is not irreducible since it
  has a $3$-cycle.\\
Right: additional skein relations for bi-webs.}
    \label{fig:sw10}
\end{figure}

The role of webs in the context of $\SL_3$ Fock-Goncharov theory
is played by \emph{bi-webs}.
These are emdedded graphs just like the usual webs \emph{\`a la}
Kuperberg, except that the boundary\linebreak[3] 
vertices are now bivariant;
interior vertices still carry a proper \hbox{$2$-coloring}. 
A bi-web is \emph{irreducible} (or \emph{non-elliptic})
if it satisfies the same conditions as before, 
and in addition avoids loops or $3$-cycles based at 
boundary points. 
We call associated $\SL(V)$-invariants \emph{bi-web invariants}. 
See Figure~\ref{fig:sw10}. 

The skein calculus for bivariant tensor diagrams employs the local
rules of Figures~\ref{fig:skein}--\ref{fig:yang-baxter}
 plus two additional relations shown in
Figure~\ref{fig:sw10}. 
The associated invariant does not change under the corresponding
transformations. 
As in Theorem~\ref{thm:confl}, the corresponding 
\emph{flattening} process is confluent:

\begin{theorem}
Any bivariant tensor diagram can be transformed by 
repeated application of skein relations 
into a linear combination of non-elliptic bi-webs. 
Furthermore, the output does not depend on the choices made.
\end{theorem}


\vspace{-.1in}

We expect Conjectures~\ref{conj:web-and-tree} and~\ref{conj:compat}
to extend to the bivariant case:

\begin{conjecture} 
\label{conj:fg}
In any $\SL_3$ Fock-Goncharov cluster algebra associated to a marked
bordered surface~$S$, 
\begin{itemize}
\vspace{-.1in}
 \item all cluster variables are bi-web invariants;
 \item a bi-web invariant is a cluster monomial if and only if it can also 
be written as a forest bivariant tensor diagram;
 \item compatibility of cluster/coefficient variables 
is equivalent to  
their product being a single bi-web invariant.
\end{itemize}
\end{conjecture}

\vspace{-.1in}

\begin{acknowledgments}
Much of this work was done at MSRI during the Fall 2012 program on
Cluster Algebras. 
We thank Greg Muller, Greg Kuperberg, and Dylan Thurston for
stimulating discussions. 

This project was partially supported by NSF grants DMS-1101152 (S.~F.)
and DMS-1068169 (P.~P.).
\end{acknowledgments}

\end{article}

\end{document}